\documentclass[oneside,reqno,english]{amsart}
\usepackage[T1]{fontenc}
\usepackage[latin9]{inputenc}
\usepackage{amstext}
\usepackage{amsthm}

\makeatletter

\providecommand{\tabularnewline}{\\}

\makeatother

\usepackage{babel}
\begin{document}
\title{Triangular Numbers Multiple of Triangular Numbers and Solutions of
Pell Equations}
\author{Vladimir Pletser}
\address{{\large{}European Space Agency (ret.)}}
\email{{\large{}Pletservladimir@gmail.com}}
\begin{abstract}
For all positive non-square integer multiplier $k$, there is an infinity
of multiples of triangular numbers which are also triangular numbers.
With a simple change of variables, these triangular numbers can be
found using solutions of Pell equations. With some conditions on parities
of fundamental solutions of the simple and generalized Pell equations,
only odd solutions of the generalized Pell equation are retained to
provide many infinitely solutions found on branches corresponding
to each of the generalized fundamental solutions. General algebraic
expressions of fundamental solutions of the Pell equations are found
for some values of the multiplier $k$ in function of the closest
natural square. Further, among the expressions of Pell equation solutions,
a set of recurrent relations is identical to those found previously
without the Pell equation solving method. It is found also that two
constants of the problem of multiples of triangular numbers are directly
related to the fundamental solutions of the simple Pell equation,
which is an unexpected result as it means that simple Pell equation
fundamental solutions in all generality, are related to constants
in recurrent relations of the problem of finding triangular numbers
multiple of other triangular numbers.
\end{abstract}

\maketitle

\section{{\large{}Introduction}}

{\large{}\label{sec1:Introduction}}{\large\par}

{\large{}Triangular numbers $T_{t}=\frac{t\left(t+1\right)}{2}$ are
figurate numbers with several interesting properties and formula (see,
e.g., \cite{Weisstein Tri Nber,Andrews}). In this paper, we investigate
triangular numbers $T_{\xi}$ that are multiples of other triangular
numbers $T_{t}$ 
\begin{equation}
T_{\xi}=kT_{t}\label{eq:1}
\end{equation}
Several authors have investigated this Diophantine equation; see,
e.g., \cite{Cunningham,Joncourt,Roegel,Dickson Sums of cubes,Chahal DSouza,Breiteig,Chahal Griffin Priddis}.
Further historical accounts can be found in \cite{Dickson Sums of cubes}.
Recently, Pletser showed \cite{Pletser Recur Rel Mult Tri Nber} that,
for non-square integer values of $k$, the four variables $t,\xi,T_{t}$
and $T_{\xi}$ can be represented by recurrent relations involving
a rank $r$ and parameters $\kappa$ and $\gamma$ which are respectively
the sum and the product of the $\left(r-1\right)^{\text{th}}$ and
the $r^{\text{th}}$ values of $t$. The rank is being defined as
the number of successive values of $t$ solutions of (\ref{eq:1})
such that their successive ratios are slowly decreasing without jumps.}{\large\par}

{\large{}We only consider solutions of (\ref{eq:1}) for $k>1$ as,
for $k=0$ and $k=1$, solutions are trivial, respectively, $\xi=0$
and $\xi=t$ for any positive integer $t$.}{\large\par}

{\large{}In this paper, we investigate how to find all solutions to
(\ref{eq:1}) using the method of resolution of the simple and generalized
Pell equations associated to (\ref{eq:1}). We show that the rank
$r$ and parameters $\kappa$ and $\gamma$ of recurrent relations
can be deduced from fundamental solutions of Pell equations. Section
2 introduces the rank $r$ and recurrent relations. Section 3 give
a short reminder on how to find solutions of Pell equations. In Section
4, Pell equation methods are applied to find all multiples of triangular
numbers that are triangular numbers. In certain cases, general expressions
of fundamental solutions of the Pell equations associated to (\ref{eq:1})
are given for values of the multiplier $k$ in function of the closest
natural square values $s^{2}$.}{\large\par}

\section{{\large{}Rank and recurrent relations}}

\label{sec2:Rank-and-recurrent}

{\large{}The Online Encyclopedia of Integer Sequences (OEIS) \cite{Sloane}
lists sequences of solutions of (\ref{eq:1}) for $k=2,3,5,6,7,8$.
Let us note first that, among all solutions, $\left(t_{0},\xi_{0}\right)=\left(0,0\right)$
is always a first solution of (\ref{eq:1}) for all non-square integer
value of $k$. }{\large\par}

{\large{}Let's consider the two cases of $k=3$ and $k=6$ yielding
the successive solution pairs as shown in Table \ref{tab1:Solutions-of-}.
We indicate also the ratios $t_{n}/t_{n-1}$ for both cases and $t_{n}/t_{n-2}$
for $k=6$. It is seen that for $k=3$, the ratio $t_{n}/t_{n-1}$
varies between close values, from 5 down to 3.737, while for $k=6$,
the ratio $t_{n}/t_{n-1}$ alternates between values 3 ... 2.385 and
4.667 ... 4.206, while the ratio $t_{n}/t_{n-2}$ decreases more regularly
from 14 to 10.029 (corresponding approximately to the product of the
alternating values of the ratio $t_{n}/t_{n-1}$). We call rank $r$
the integer value such that $t_{n}/t_{n-r}$ is approximately constant
or, better, decreases regularly without jumps (a more precise definition
is given further). So, here, the case $k=3$ has rank $r=1$ and the
case $k=6$ has rank $r=2$.}{\large\par}

{\large{}}
\begin{table}
{\large{}\caption{\label{tab1:Solutions-of-}Solutions of (\ref{eq:1}) for $k=3,6$}
}{\large\par}
\centering{}{\large{}}%
\begin{tabular}{|c|rrl|rrll|}
\cline{2-8} \cline{3-8} \cline{4-8} \cline{5-8} \cline{6-8} \cline{7-8} \cline{8-8} 
\multicolumn{1}{c|}{} &  & {\large{}$k=3$} &  &  &  & {\large{}$k=6$} & \tabularnewline
\hline 
{\large{}$n$} & {\large{}$t_{n}$} & {\large{}$\xi_{n}$} & {\large{}$t_{n}/t_{n-1}$} & {\large{}$t_{n}$} & {\large{}$\xi_{n}$} & {\large{}$t_{n}/t_{n-1}$} & {\large{}$t_{n}/t_{n-2}$}\tabularnewline
\hline 
\hline 
{\large{}0} & {\large{}0} & {\large{}0} &  & {\large{}0} & {\large{}0} &  & \tabularnewline
\hline 
{\large{}1} & {\large{}1} & {\large{}2} & {\large{}--} & {\large{}1} & {\large{}3} & {\large{}--} & {\large{}--}\tabularnewline
\hline 
{\large{}2} & {\large{}5} & {\large{}9} & {\large{}5} & {\large{}3} & {\large{}8} & {\large{}3} & {\large{}--}\tabularnewline
\hline 
{\large{}3} & {\large{}20} & {\large{}35} & {\large{}4} & {\large{}14} & {\large{}35} & {\large{}4.66667} & {\large{}14}\tabularnewline
\hline 
{\large{}4} & {\large{}76} & {\large{}132} & {\large{}3.8} & {\large{}34} & {\large{}84} & {\large{}2.42857} & {\large{}11.33333}\tabularnewline
\hline 
{\large{}5} & {\large{}285} & {\large{}494} & {\large{}3.75} & {\large{}143} & {\large{}351} & {\large{}4.20588} & {\large{}10.21429}\tabularnewline
\hline 
{\large{}6} & {\large{}1065} & {\large{}1845} & {\large{}3.73684} & {\large{}341} & {\large{}836} & {\large{}2.38461} & {\large{}10.02941}\tabularnewline
\hline 
\end{tabular}{\large\par}
\end{table}
{\large{}Pletser showed \cite{Pletser Recur Rel Mult Tri Nber} that
the rank $r$ is the index of $t_{r}$ and $\xi_{r}$ solutions of
(\ref{eq:1}) such that}{\large\par}

{\large{}
\begin{equation}
\kappa=t_{r}+t_{r-1}=\xi_{r}-\xi_{r-1}-1\label{eq:3.2}
\end{equation}
The rank $r$ is also such that the ratio $t_{2r}/t_{r}$, corrected
by the ratio $t_{r-1}/t_{r}$, is equal to a constant $2\kappa+3$
\begin{equation}
\frac{t_{2r}-t_{r-1}}{t_{r}}=2\kappa+3\label{eq:3-0}
\end{equation}
For example, for $k=6$ and $r=2$, $\kappa=t_{2}+t_{1}=3+1=4$, and
$\kappa=\xi_{2}-\xi_{1}-1=8-3-1=4$, yielding $2\kappa+3=11$.}{\large\par}

{\large{}Pletser found \cite{Pletser Recur Rel Mult Tri Nber} four
recurrent equations for $t_{n},\xi_{n},T_{t_{n}}$ and $T_{\xi_{n}}$for
each non-square integer value of $k$ }{\large\par}

{\large{}
\begin{align}
t_{n} & =2\left(\kappa+1\right)t_{n-r}-t_{n-2r}+\kappa\label{eq:3.3}\\
\xi_{n} & =2\left(\kappa+1\right)\xi_{n-r}-\xi_{n-2r}+\kappa\label{eq:3.3-1}\\
T_{t_{n}} & =\left(4\left(\kappa+1\right)^{2}-2\right)T_{t_{n-r}}-T_{t_{n-2r}}+\left(T_{\kappa}-\gamma\right)\label{eq:3.3-2}\\
T_{\xi_{n}} & =\left(4\left(\kappa+1\right)^{2}-2\right)T_{\xi_{n-r}}-T_{\xi_{n-2r}}+k\left(T_{\kappa}-\gamma\right)\label{eq:3.3-3}
\end{align}
where coefficients are functions of two constants $\kappa$ and $\gamma$,
respectively the sum (\ref{eq:3.2}) and the product $\gamma=t_{r-1}t_{r}$.
Note that these four relations are independent from the value of $k$. }{\large\par}

\section{{\large{}Pell equations: A Reminder}}

\label{sec3:Pell-equations:-A}

\noindent {\large{}The Diophantine bivariate quadratic equation}{\large\par}

{\large{}
\begin{equation}
X^{2}-DY^{2}=N,\label{eq:2}
\end{equation}
}{\large\par}

\noindent {\large{}with integers $X,Y,D,N$ and square free $D$,
is called the Pell equation. Several mathematicians have investigated
this equation (see historical accounts in \cite{Weil,Lenstra,Lemmermeyer,O'Connor Robertson,Dickson Pell eq}),
Treatments and solutions are described in several classical text books
(see e.g. \cite{Nagell,Jacobson Williams,Weisstein Pell eq,Robertson}
and references therein). We remind here some general formulas and
how to calculate solutions. Details can be found in references.}{\large\par}

\noindent {\large{}For $N=1$, (\ref{eq:2}) is called the simple
Pell equation 
\begin{equation}
x^{2}-Dy^{2}=1\label{eq:3}
\end{equation}
This equation admits the obvious trivial solution $\left(x_{0},y_{0}\right)=\left(1,0\right)$
and infinitely many solutions given by
\begin{align}
\left(x_{n},y_{n}\right)= & \left(\frac{\left(x_{f}+\sqrt{D}y_{f}\right)^{n}+\left(x_{f}-\sqrt{D}y_{f}\right)^{n}}{2},\right.\nonumber \\
 & \left.\frac{\left(x_{f}+\sqrt{D}y_{f}\right)^{n}-\left(x_{f}-\sqrt{D}y_{f}\right)^{n}}{2\sqrt{D}}\right)\label{eq:4}
\end{align}
where $n$ are positive integers and $\left(x_{f},y_{f}\right)$ is
the least solution to (\ref{eq:3}), i.e. the smallest integer solution
different from the trivial solution, $x_{f}>1,y_{f}>0$. We call this
least solution the fundamental solution. Obviously, having found the
fundamental solution $\left(x_{f},y_{f}\right)$ yields directly three
other solutions, $\left(-x_{f},y_{f}\right),\left(x_{f},-y_{f}\right),\left(-x_{f},-y_{f}\right)$.}{\large\par}

\noindent {\large{}Lagrange devised a method to find the fundamental
solution, based on the continued fraction expansion of the quadratic
irrational $\sqrt{D}$, that can be summarized as follows. One computes
the $j^{\textnormal{th}}$ convergent $\left(p_{j}/q_{j}\right)$
of the continued fraction $\left[\alpha_{0};\alpha_{1},\ldots,\alpha_{j},\alpha_{j+1},\ldots\right]$
of $\sqrt{D}$, with $\alpha_{0}=\left\lfloor \sqrt{D}\right\rfloor $,
i.e., the greatest integer $\leq\sqrt{D}$. This continued fraction
becomes periodic after the following term, $\alpha_{j+1}=2\alpha_{0}$
if $\sqrt{D}$ is a quadratic irrational. The recurrence relations
\[
p_{i}=\alpha_{i}p_{i-1}+p_{i-2}\,\,,\,\,\,q_{i}=\alpha_{i}q_{i-1}+q_{i-2}
\]
yield the terms $p_{i}$ and $q_{i}$ of the convergent, with $p_{-2}=0,p_{-1}=1,q_{-2}=1,q_{-1}=0$.
The fundamental solution is then $\left(x_{f},y_{f}\right)=\left(p_{j},q_{j}\right)$
if $j$ is odd, or $\left(x_{f},y_{f}\right)=\left(p_{2j+1},q_{2j+1}\right)$
if $j$ is even. }{\large\par}

\noindent {\large{}For $N\neq1$, (\ref{eq:2}) is called the generalized
Pell equation, which can have either no solution, or one, or several
fundamental solutions $\left(X_{f_{i}},Y_{f_{i}}\right)$, with positive
integers $i$ such that $1\leq i\leq\rho$, where $\rho$ is the total
number of fundamental solutions admitted by (\ref{eq:2}). All integer
solutions, if they exist, are found on double infinite branches that
can be expressed in terms of the fundamental solution(s) $\left(X_{f_{i}},Y_{f_{i}}\right)$
and $\left(-X_{f_{i}},Y_{f_{i}}\right)$. Methods to calculate the
fundamental solution(s) of the generalized Pell equation (see e.g.
\cite{Nagell,Jacobson Williams,Robertson,Lagrange,Chrystal,Mollin,Matthew Integers,Matthews BCMATH}
and references therein) are all based on Lagrange's method of continued
fractions, sometime adapted (see e.g. \cite{Pletser Cont'd Fract JNT}).
The nearest integer continued fraction method and the Lagrange-Mollin-Matthews
method \cite{Matthews BCMATH} are used further to calculate the fundamental
solutions of respectively, the simple and the generalized Pell equations. }{\large\par}

{\large{}Once fundamental solutions are known, the other solutions
$\left(X_{n},Y_{n}\right)$ of (\ref{eq:2}) are calculated by 
\begin{equation}
X_{n}+\sqrt{D}Y_{n}=\pm\left(X_{f_{i}}+\sqrt{D}Y_{f_{i}}\right)\left(x_{f}+\sqrt{D}y_{f}\right)^{n}\label{eq:5}
\end{equation}
for a proper choice of sign $\pm$ \cite{Robertson}, yielding respectively,
for $n=0,1,2$ (assuming a $+$ sign),
\begin{align}
\left(X_{0},Y_{0}\right)= & \left(X_{f_{i}},Y_{f_{i}}\right)\label{eq:6}\\
\left(X_{1},Y_{1}\right)= & \left(X_{f_{i}}x_{f}+DY_{f_{i}}y_{f},X_{f_{i}}y_{f}+Y_{f_{i}}x_{f}\right)\label{eq:7}\\
\left(X_{2},Y_{2}\right)= & \left(X_{f_{i}}\left(x_{f}^{2}+Dy_{f}^{2}\right)+2DY_{f_{i}}x_{f}y_{f},\right.\nonumber \\
 & \left.Y_{f_{i}}\left(x_{f}^{2}+Dy_{f}^{2}\right)+2X_{f_{i}}x_{f}y_{f}\right)\label{eq:8}
\end{align}
Note that, for each value of $n$, one can have several (up to $\rho$)
solutions depending on the different values of the generalized fundamental
solutions $\left(X_{f_{i}},Y_{f_{i}}\right)$.}{\large\par}

{\large{}The other solutions $\left(X_{n},Y_{n}\right)$ of (\ref{eq:2})
can also be represented by recurrence relations
\begin{equation}
\left(X_{n},Y_{n}\right)=\left(x_{f}X_{n-1}+Dy_{f}Y_{n-1},x_{f}Y_{n-1}+y_{f}X_{n-1}\right)\label{eq:9}
\end{equation}
that can also be written as
\begin{equation}
\left(X_{n},Y_{n}\right)=\left(2x_{f}X_{n-1}-X_{n-2},2x_{f}Y_{n-1}-Y_{n-2}\right)\label{eq:10}
\end{equation}
}{\large\par}

\noindent {\large{}or by Chebyshev polynomials of the first kind $T_{n-1}\left(x_{f}\right)$
and of the second kind $U_{n-2}\left(x_{f}\right)$, evaluated at
$x_{f}$ (see \cite{Pletser Gen Pell eq Chebysh Polyn}),
\begin{align}
\left(X_{n},Y_{n}\right)= & \left(X_{f_{i}}T_{n-1}\left(x_{f}\right)+DY_{f_{i}}y_{f}U_{n-2}\left(x_{f}\right),\right.\nonumber \\
 & \left.X_{f_{i}}y_{f}U_{n-2}\left(x_{f}\right)+Y_{f_{i}}T_{n-1}\left(x_{f}\right)\right)\label{eq:11}
\end{align}
One notices that the second recurrent relations (\ref{eq:10}) is
similar in form to recurrent relations (\ref{eq:3.3}) and (\ref{eq:3.3-1})
found in \cite{Pletser Recur Rel Mult Tri Nber}.}{\large\par}

\section{{\large{}Pell Equations and Multiples of Triangular Numbers}}

\label{sec4:Pell-Equations-and}

\subsection{Solutions of Pell equations}

\label{subsec4.1:Solutions-of-Pell}

{\large{}For non-square integers $k$ and with the change of variables
\begin{equation}
\left(X,Y\right)=\left(2\xi+1,2t+1\right)\label{eq:12}
\end{equation}
(\ref{eq:1}) becomes a generalized Pell equation \cite{Chahal DSouza,Breiteig},
with $D=k$ and $N=1-k$ negative, as $k>1$,
\begin{equation}
X^{2}-kY^{2}=1-k\label{eq:13}
\end{equation}
and the associated simple Pell equation reads} {\large{}
\begin{equation}
x^{2}-ky^{2}=1\label{eq:14}
\end{equation}
Odd solutions $\left(X,Y\right)$ of (\ref{eq:13}) provide then pairs
$\left(\xi,t\right)$, solutions of (\ref{eq:1}). Following the procedure
of Section 2, the fundamental solutions of the simple and generalized
Pell equations are calculated and shown in Tables \ref{tab2 Fund Sol Pell Eq}
to (\ref{tab4:Fundamental-solutions-of}) for non-square $k$ between
2 and 102. The second and third columns give the rank $r$ found in
\cite{Pletser Recur Rel Mult Tri Nber} and the total number $\rho$
of fundamental solutions of the generalized Pell equation.The fourth
column shows the single fundamental solution of the simple Pell equations;
the fifth and sixth columns give the fundamental solutions of the
generalized Pell equations, the fifth column for those solutions with
both $X_{f_{i}}$ and $Y_{f_{i}}$ odd or having different parities,
while the sixth column give those solutions with both $X_{f_{i}}$
and $Y_{f_{i}}$ even (except for $k=56$, see discussion further).}{\large\par}

{\large{}}
\begin{table}
{\large{}\caption{\label{tab2 Fund Sol Pell Eq}Fundamental solutions of simple (\ref{eq:13})
and generalized (\ref{eq:14}) Pell equations}
}{\large\par}
\centering{}{\large{}}%
\begin{tabular}{|r|r|r|l|l|l|}
\hline 
$k$ & $r$ & $\rho$ & $\left(x_{f},y_{f}\right)$ & $\left(X_{f_{i}},Y_{f_{i}}\right)$ & \tabularnewline
\hline 
\hline 
2 & 1 & 1 & $\left(3,2\right)$ & $\left(1,1\right)$ & \tabularnewline
\hline 
3 & 1 & 1 & $\left(2,1\right)$ & $\left(1,1\right)$ & \tabularnewline
\hline 
5 & 2 & 3 & $\left(9,4\right)$ & $\left(\pm1,1\right)$ & $\left(4,2\right)$\tabularnewline
\hline 
6 & 2 & 2 & $\left(5,2\right)$ & $\left(\pm1,1\right)$ & \tabularnewline
\hline 
7 & 2 & 2 & $\left(8,3\right)$ & $\left(\pm1,1\right)$ & \tabularnewline
\hline 
8 & 2 & 2 & $\left(3,1\right)$ & $\left(\pm1,1\right)$ & \tabularnewline
\hline 
10 & 3 & 3 & $\left(19,6\right)$ & $\left(\pm1,1\right),\left(9,3\right)$ & \tabularnewline
\hline 
11 & 2 & 2 & $\left(10,3\right)$ & $\left(\pm1,1\right)$ & \tabularnewline
\hline 
12 & 2 & 2 & $\left(7,2\right)$ & $\left(\pm1,1\right)$ & \tabularnewline
\hline 
13 & 4 & 6 & $\left(649,180\right)$ & $\left(\pm1,1\right),\left(\pm25,7\right)$ & $\left(\pm14,4\right)$\tabularnewline
\hline 
14 & 2 & 2 & $\left(15,4\right)$ & $\left(\pm1,1\right)$ & \tabularnewline
\hline 
15 & 2 & 2 & $\left(4,1\right)$ & $\left(\pm1,1\right)$ & \tabularnewline
\hline 
17 & 2 & 3 & $\left(33,8\right)$ & $\left(\pm1,1\right)$ & $\left(16,4\right)$\tabularnewline
\hline 
18 & 2 & 2 & $\left(17,4\right)$ & $\left(\pm1,1\right)$ & \tabularnewline
\hline 
19 & 3 & 3 & $\left(170,39\right)$ & $\left(\pm1,1\right),\left(39,9\right)$ & \tabularnewline
\hline 
20 & 2 & 2 & $\left(9,2\right)$ & $\left(\pm1,1\right)$ & \tabularnewline
\hline 
21 & 4 & 6 & $\left(55,12\right)$ & $\left(\pm1,1\right),\left(\pm13,3\right)$ & $\left(\pm8,2\right)$\tabularnewline
\hline 
22 & 4 & 4 & $\left(197,42\right)$ & $\left(\pm1,1\right),\left(\pm23,5\right)$ & \tabularnewline
\hline 
23 & 2 & 2 & $\left(24,5\right)$ & $\left(\pm1,1\right)$ & \tabularnewline
\hline 
24 & 2 & 2 & $\left(5,1\right)$ & $\left(\pm1,1\right)$ & \tabularnewline
\hline 
26 & 3 & 3 & $\left(51,10\right)$ & $\left(\pm1,1\right),\left(25,5\right)$ & \tabularnewline
\hline 
27 & 2 & 2 & $\left(26,5\right)$ & $\left(\pm1,1\right)$ & \tabularnewline
\hline 
28 & 4 & 4 & $\left(127,24\right)$ & $\left(\pm1,1\right),\left(\pm15,3\right)$ & \tabularnewline
\hline 
29 & 4 & 6 & $\left(9801,1820\right)$ & $\left(\pm1,1\right),\left(\pm59,11\right)$ & $\left(\pm86,16\right)$\tabularnewline
\hline 
30 & 2 & 2 & $\left(11,2\right)$ & $\left(\pm1,1\right)$ & \tabularnewline
\hline 
31 & 4 & 4 & $\left(1520,273\right)$ & $\left(\pm1,1\right),\left(\pm61,11\right)$ & \tabularnewline
\hline 
32 & 2 & 2 & $\left(17,3\right)$ & $\left(\pm1,1\right)$ & \tabularnewline
\hline 
33 & 2 & 4 & $\left(23,4\right)$ & $\left(\pm1,1\right)$ & $\left(\pm10,2\right)$\tabularnewline
\hline 
34 & 2 & 2 & $\left(35,6\right)$ & $\left(\pm1,1\right)$ & \tabularnewline
\hline 
35 & 2 & 2 & $\left(6,1\right)$ & $\left(\pm1,1\right)$ & \tabularnewline
\hline 
37 & 2 & 3 & $\left(73,12\right)$ & $\left(\pm1,1\right)$ & $\left(36,6\right)$\tabularnewline
\hline 
38 & 2 & 2 & $\left(37,6\right)$ & $\left(\pm1,1\right)$ & \tabularnewline
\hline 
39 & 2 & 2 & $\left(25,4\right)$ & $\left(\pm1,1\right)$ & \tabularnewline
\hline 
40 & 4 & 4 & $\left(19,3\right)$ & $\left(\pm1,1\right),\left(\pm11,2\right)$ & \tabularnewline
\hline 
41 & 4 & 4 & $\left(2049,320\right)$ & $\left(\pm1,1\right),\left(\pm83,13\right)$ & \tabularnewline
\hline 
42 & 2 & 2 & $\left(13,2\right)$ & $\left(\pm1,1\right)$ & \tabularnewline
\hline 
43 & 4 & 4 & $\left(3482,531\right)$ & $\left(\pm1,1\right),\left(\pm85,13\right)$ & \tabularnewline
\hline 
44 & 2 & 2 & $\left(199,30\right)$ & $\left(\pm1,1\right)$ & \tabularnewline
\hline 
45 & 4 & 6 & $\left(161,24\right)$ & $\left(\pm1,1\right),\left(\pm19,3\right)$ & $\left(\pm26,4\right)$\tabularnewline
\hline 
\end{tabular}{\large\par}
\end{table}
{\large\par}

\begin{table}
\caption{\label{tab3:Fundamental-solutions-of-1}Fundamental solutions of simple
(\ref{eq:13}) and generalized (\ref{eq:14}) Pell equations}

\centering{}{\large{}}%
\begin{tabular}{|r|r|r|l|l|l|}
\hline 
$k$ & $r$ & $\rho$ & $\left(x_{f},y_{f}\right)$ & $\left(X_{f_{i}},Y_{f_{i}}\right)$ & \tabularnewline
\hline 
\hline 
46 & 6 & 6 & $\left(24335,3588\right)$ & $\left(\pm1,1\right),\left(\pm47,7\right),\left(\pm183,27\right)$ & \tabularnewline
\hline 
47 & 2 & 2 & $\left(48,7\right)$ & $\left(\pm1,1\right)$ & \tabularnewline
\hline 
48 & 2 & 2 & $\left(7,1\right)$ & $\left(\pm1,1\right)$ & \tabularnewline
\hline 
50 & 3 & 3 & $\left(99,14\right)$ & $\left(\pm1,1\right),\left(49,7\right)$ & \tabularnewline
\hline 
51 & 3 & 3 & $\left(50,7\right)$ & $\left(\pm1,1\right),\left(35,5\right)$ & \tabularnewline
\hline 
52 & 4 & 4 & $\left(649,90\right)$ & $\left(\pm1,1\right),\left(\pm79,11\right)$ & \tabularnewline
\hline 
53 & 4 & 6 & $\left(66249,9100\right)$ & $\left(\pm1,1\right),\left(\pm211,29\right)$ & $\left(\pm160,22\right)$\tabularnewline
\hline 
54 & 2 & 2 & $\left(485,66\right)$ & $\left(\pm1,1\right)$ & \tabularnewline
\hline 
55 & 4 & 4 & $\left(89,12\right)$ & $\left(\pm1,1\right),\left(\pm21,3\right)$ & \tabularnewline
\hline 
56 & 2 & 4 & $\left(15,2\right)$ & $\left(\pm1,1\right)$ & $\left(\pm13,2\right)$\tabularnewline
\hline 
57 & 4 & 4 & $\left(151,20\right)$ & $\left(\pm1,1\right),\left(\pm37,5\right)$ & \tabularnewline
\hline 
58 & 4 & 4 & $\left(19603,2574\right)$ & $\left(\pm1,1\right),\left(\pm175,23\right)$ & \tabularnewline
\hline 
59 & 2 & 2 & $\left(530,69\right)$ & $\left(\pm1,1\right)$ & \tabularnewline
\hline 
60 & 2 & 2 & $\left(31,4\right)$ & $\left(\pm1,1\right)$ & \tabularnewline
\hline 
61 & 8 & 12 & $\left(1766319049,\right.$ & $\left(\pm1,1\right),\left(\pm367,47\right),$ & $\left(\pm62,8\right),$\tabularnewline
 &  &  & $\left.226153980\right)$ & $\left(\pm6709,859\right),\left(\pm94793,12137\right)$ & $\left(\pm5186,664\right)$\tabularnewline
\hline 
62 & 2 & 2 & $\left(63,8\right)$ & $\left(\pm1,1\right)$ & \tabularnewline
\hline 
63 & 2 & 2 & $\left(8,1\right)$ & $\left(\pm1,1\right)$ & \tabularnewline
\hline 
65 & 2 & 5 & $\left(129,16\right)$ & $\left(\pm1,1\right)$ & $\left(\pm14,2\right),\left(64,8\right)$\tabularnewline
\hline 
66 & 4 & 4 & $\left(65,8\right)$ & $\left(\pm1,1\right),\left(\pm23,3\right)$ & \tabularnewline
\hline 
67 & 4 & 4 & $\left(48842,5967\right)$ & $\left(\pm1,1\right),\left(\pm401,49\right)$ & \tabularnewline
\hline 
68 & 2 & 2 & $\left(33,4\right)$ & $\left(\pm1,1\right)$ & \tabularnewline
\hline 
69 & 4 & 6 & $\left(7775,936\right)$ & $\left(\pm1,1\right),\left(\pm91,11\right)$ & $\left(\pm116,14\right)$\tabularnewline
\hline 
70 & 4 & 4 & $\left(251,30\right)$ & $\left(\pm1,1\right),\left(\pm41,5\right)$ & \tabularnewline
\hline 
71 & 4 & 4 & $\left(3480,413\right)$ & $\left(\pm1,1\right),\left(\pm143,17\right)$ & \tabularnewline
\hline 
72 & 2 & 2 & $\left(17,2\right)$ & $\left(\pm1,1\right)$ & \tabularnewline
\hline 
73 & 6 & 6 & $\left(2281249,267000\right)$ & $\left(\pm1,1\right),\left(\pm145,17\right),\left(\pm1461,171\right)$ & \tabularnewline
\hline 
74 & 2 & 2 & $\left(3699,430\right)$ & $\left(\pm1,1\right)$ & \tabularnewline
\hline 
75 & 2 & 2 & $\left(26,3\right)$ & $\left(\pm1,1\right)$ & \tabularnewline
\hline 
76 & 6 & 6 & $\left(57799,6630\right)$ & $\left(\pm1,1\right),\left(\pm113,13\right),\left(\pm305,35\right)$ & \tabularnewline
\hline 
77 & 4 & 6 & $\left(351,40\right)$ & $\left(\pm1,1\right),\left(\pm43,5\right)$ & $\left(\pm34,4\right)$\tabularnewline
\hline 
78 & 4 & 4 & $\left(53,6\right)$ & $\left(\pm1,1\right),\left(\pm25,3\right)$ & \tabularnewline
\hline 
79 & 2 & 2 & $\left(80,9\right)$ & $\left(\pm1,1\right)$ & \tabularnewline
\hline 
80 & 2 & 2 & $\left(9,1\right)$ & $\left(\pm1,1\right)$ & \tabularnewline
\hline 
82 & 3 & 3 & $\left(163,18\right)$ & $\left(\pm1,1\right),\left(81,9\right)$ & \tabularnewline
\hline 
83 & 2 & 2 & $\left(82,9\right)$ & $\left(\pm1,1\right)$ & \tabularnewline
\hline 
84 & 2 & 2 & $\left(55,6\right)$ & $\left(\pm1,1\right)$ & \tabularnewline
\hline 
85 & 8 & 12 & $\left(285769,\right.$ & $\left(\pm1,1\right),\left(\pm101,11\right),$ & $\left(\pm16,2\right),$\tabularnewline
 &  &  & $\left.30996\right)$ & $\left(\pm341,37\right),\left(\pm1429,155\right)$ & $\left(\pm424,46\right)$\tabularnewline
\hline 
\end{tabular}{\large\par}
\end{table}

\begin{table}
\caption{\label{tab4:Fundamental-solutions-of}Fundamental solutions of simple
(\ref{eq:13}) and generalized (\ref{eq:14}) Pell equations}

\centering{}{\large{}}%
\begin{tabular}{|r|r|r|l|l|l|}
\hline 
$k$ & $r$ & $\rho$ & $\left(x_{f},y_{f}\right)$ & $\left(X_{f_{i}},Y_{f_{i}}\right)$ & \tabularnewline
\hline 
\hline 
86 & 4 & 4 & $\left(10405,1122\right)$ & $\left(\pm1,1\right),\left(\pm343,37\right)$ & \tabularnewline
\hline 
87 & 2 & 2 & $\left(28,3\right)$ & $\left(\pm1,1\right)$ & \tabularnewline
\hline 
88 & 4 & 4 & $\left(197,21\right)$ & $\left(\pm1,1\right),\left(\pm65,7\right)$ & \tabularnewline
\hline 
89 & 4 & 4 & $\left(500001,53000\right)$ & $\left(\pm1,1\right),\left(\pm179,19\right)$ & \tabularnewline
\hline 
90 & 2 & 2 & $\left(19,2\right)$ & $\left(\pm1,1\right)$ & \tabularnewline
\hline 
91 & 6 & 6 & $\left(1574,165\right)$ & $\left(\pm1,1\right),\left(\pm27,3\right),\left(\pm181,19\right)$ & \tabularnewline
\hline 
92 & 4 & 4 & $\left(1151,120\right)$ & $\left(\pm1,1\right),\left(\pm47,5\right)$ & \tabularnewline
\hline 
93 & 4 & 6 & $\left(12151,1260\right)$ & $\left(\pm1,1\right),\left(\pm125,13\right)$ & $\left(\pm154,16\right)$\tabularnewline
\hline 
94 & 4 & 4 & $\left(2143295,\right.$ & $\left(\pm1,1\right),\left(\pm281,29\right)$ & \tabularnewline
 &  &  & $\left.221064\right)$ &  & \tabularnewline
\hline 
95 & 2 & 2 & $\left(39,4\right)$ & $\left(\pm1,1\right)$ & \tabularnewline
\hline 
96 & 4 & 4 & $\left(49,5\right)$ & $\left(\pm1,1\right),\left(\pm17,2\right)$ & \tabularnewline
\hline 
97 & 4 & 8 & $\left(62809633,\right.$ & $\left(\pm1,1\right),\left(\pm581,59\right)$ & $\left(\pm98,10\right),$\tabularnewline
 &  &  & $\left.6377352\right)$ &  & $\left(\pm12902,1310\right)$\tabularnewline
\hline 
98 & 2 & 2 & $\left(99,10\right)$ & $\left(\pm1,1\right)$ & \tabularnewline
\hline 
99 & 2 & 2 & $\left(10,1\right)$ & $\left(\pm1,1\right)$ & \tabularnewline
\hline 
101 & 2 & 3 & $\left(201,20\right)$ & $\left(\pm1,1\right)$ & $\left(100,10\right)$\tabularnewline
\hline 
102 & 2 & 2 & $\left(101,10\right)$ & $\left(\pm1,1\right)$ & \tabularnewline
\hline 
\end{tabular}{\large\par}
\end{table}

{\large{}From these Tables, we deduce the following.}{\large\par}

{\large{}First, the rank of solutions of (\ref{eq:1}) is equal to,
or less than, the total number of fundamental solutions of the generalized
Pell equations, $r\leq\rho$, as was expected.}{\large\par}

{\large{}Second, for all the single fundamental solutions $\left(x_{f},y_{f}\right)$
of the simple Pell equation, both $x_{f}$ and $y_{f}$ are of different
parities, i.e., one is odd, the other even (except for some cases
of $k\equiv0\left(\text{mod}8\right)$, where both $x_{f}$ and $y_{f}$
are odd; see further). It is easy to see why: for (\ref{eq:14}) to
hold, the following three conditions must hold:}{\large\par}

{\large{}(C1) $x_{f}$ and $y_{f}$ can not be simultaneously even,
whatever the value of $k$ is; }{\large\par}

{\large{}(C2) if $k$ is even, $x_{f}$ must necessarily be odd and
$y_{f}$ can be either even or odd; }{\large\par}

{\large{}(C3) if $k$ is odd, $x_{f}$ and $y_{f}$ must have different
parities, one odd and the other even.}{\large\par}

{\large{}Third, the sets of fundamental solutions of the generalized
Pell equation always include the two fundamental solutions $\left(X_{f_{1}},Y_{f_{1}}\right)=\left(1,1\right)$
and $\left(X_{f_{2}},Y_{f_{2}}\right)=\left(-1,1\right)$, which is
quite obvious from (\ref{eq:13}). The only two exceptions are for
the cases $k=2$ and $3$. Although $\left(-1,1\right)$ is also a
solution to (\ref{eq:13}) for these two cases, it does not bring
a new branch of solutions calculated by (\ref{eq:5}) to (\ref{eq:7})
different from the one obtained with $\left(1,1\right)$. Therefore,
there is only one fundamental solution, i.e., $\rho=1$ for these
two cases. Furthermore, the two pairs $\left(1,-1\right)$ and $\left(-1,-1\right)$
are also solutions of (\ref{eq:13}), but they do not yield new branches
of solutions different from those obtained with $\left(-1,1\right)$
and $\left(1,1\right)$. }{\large\par}

{\large{}Fourth, all generalized fundamental solutions $\left(X_{f_{i}},Y_{f_{i}}\right)$
with $i>2$, i.e., other than $\left(\pm1,1\right)$, have both $X_{f_{i}}$
and $Y_{f_{i}}$ odd, except for $k=40,96,208,\ldots$ where $Y_{f_{i}}$
is even.}{\large\par}

{\large{}Fifth, the generalized fundamental solutions with both $X_{f_{i}}$
and $Y_{f_{i}}$ even are shown separately as they do not bring any
solutions to (\ref{eq:1}), and there are $\rho-r$ such solutions.}{\large\par}

{\large{}With the two generalized fundamental solutions $\left(X_{f_{1,2}},Y_{f_{1,2}}\right)=\left(\pm1,1\right)$,
one has from (\ref{eq:6}) $\left(X_{0_{1,2}},Y_{0_{1,2}}\right)=\left(\pm1,1\right)$
and it yields the two trivial solutions $\left(\xi_{0_{1,2}},t_{0_{1,2}}\right)=\left(\frac{\pm1-1}{2},\frac{1-1}{2}\right)=\left(0,0\right)$
and $\left(-1,0\right)$ of (\ref{eq:1}). The next generalized solution
(\ref{eq:7}) reads
\[
\left(X_{1_{1,2}},Y_{1_{1,2}}\right)=\left(\left(\pm x_{f}+ky_{f}\right),\left(\pm y_{f}+x_{f}\right)\right)
\]
yielding, from (\ref{eq:12}),
\begin{equation}
\left(\xi_{1,2},t_{1,2}\right)=\left(\frac{\pm x_{f}+ky_{f}-1}{2},\frac{\pm y_{f}+x_{f}-1}{2}\right)\label{eq:15}
\end{equation}
 with both terms integers under the three conditions C1 to C3 above. }{\large\par}

{\large{}For other generalized fundamental solutions $\left(X_{f_{i}},Y_{f_{i}}\right)$
(with $i>2$) different from $\left(\pm1,1\right)$, one has from
(\ref{eq:6}) $\left(X_{0_{i}},Y_{0_{i}}\right)=\left(X_{f_{i}},Y_{f_{i}}\right)$,
yielding
\begin{equation}
\left(\xi_{0_{1,2}},t_{0_{1,2}}\right)=\left(\frac{X_{f_{i}}-1}{2},\frac{Y_{f_{i}}-1}{2}\right)\label{eq:15-1}
\end{equation}
integer solutions of (\ref{eq:1}) if $X_{f_{i}}$ and $Y_{f_{i}}$
are both odd. The next generalized solution (\ref{eq:7}) reads $\left(X_{1_{i}},Y_{1_{i}}\right)=\left(X_{f_{i}}x_{f}+kY_{f_{i}}y_{f},X_{f_{i}}y_{f}+Y_{f_{i}}x_{f}\right)$,
yielding
\begin{equation}
\left(\xi_{1_{i}},t_{1_{i}}\right)=\left(\frac{X_{f_{i}}x_{f}+kY_{f_{i}}y_{f}-1}{2},\frac{X_{f_{i}}y_{f}+Y_{f_{i}}x_{f}-1}{2}\right)\label{eq:16}
\end{equation}
 One sees clearly that $X_{f_{i}}$ and $Y_{f_{i}}$ can not be simultaneously
even for $\xi_{1_{i}}$ and $t_{1_{i}}$ to be integers. For $X_{f_{i}}$
and $Y_{f_{i}}$ both odd, the three conditions C1 to C3 above ensure
that $\xi_{1_{i}}$ and $t_{1_{i}}$ are integers.}{\large\par}

{\large{}For the cases of $X_{f_{i}}$ odd and $Y_{f_{i}}$ even,
like for $k=40$ and $96$ in Tables \ref{tab2 Fund Sol Pell Eq}
to \ref{tab4:Fundamental-solutions-of}, one has that $x_{f}$ and
$y_{f}$ must be simultaneously odd and, by condition C2 above, $k$
must be even for (\ref{eq:16}) to provide integer solutions.}{\large\par}

{\large{}Finally, for all single fundamental solutions $\left(x_{f},y_{f}\right)$
of the simple Pell equation with both $x_{f}$ and $y_{f}$ odd,they
appear for most of the values of $k$ such that $k\equiv0\left(\text{mod}8\right)$.
The exceptions to this are for $k=56$, 72, 112, 184, 240, 248, 264,
272, 376, ..., i.e., for some values of $k$ such that $k\equiv\pm8,\pm16\left(\text{mod}64\right)$
(but not all), where $y_{f}$ is even. In these cases, one has that
$k$ and $y_{f}$ are both even, then $Y_{f_{i}}$ can not be even
for (\ref{eq:16}) to provide integer solutions. If this is not the
case, i.e., if $Y_{f_{i}}$ is even, then the generalized fundamental
solutions $\left(X_{f_{i}},Y_{f_{i}}\right)$ must be discarded as
it does not provide integer solutions for $t$ in (\ref{eq:16}).}{\large\par}

{\large{}For the general case of $k\equiv0\left(\text{mod}8\right)$,
the fact that $y_{f}$ is not odd can be explained as follows. As
$k\equiv0\left(\text{mod}8\right)$ is not square free, the simple
Pell equation (\ref{eq:14}) can be simplified posing $k=c^{2}k^{\prime}$,
with $k^{\prime}$ square free, yielding
\begin{equation}
x^{2}-k^{\prime}y^{\prime}{}^{2}=1\label{eq:17}
\end{equation}
with $y^{\prime}=cy$. The fundamental solution $\left(x_{f},y_{f}^{\prime}\right)$
of (\ref{eq:17}) yields then the fundamental solution $\left(x_{f},y_{f}\right)=\left(x_{f},\frac{y_{f}^{\prime}}{c}\right)$
of (\ref{eq:14}). For example, for $k=8$, let $k^{\prime}=2$ and
$c=2,$(\ref{eq:17}) yields $\left(x_{f},y_{f}^{\prime}\right)=\left(3,2\right)$
and $\left(x_{f},\frac{y_{f}^{\prime}}{c}\right)=\left(x_{f},y_{f}\right)=\left(3,1\right)$.
For most of the cases of $k$ such that $k\equiv0\left(\text{mod}8\right)$,
$y_{f}^{\prime}$ is divisible by c such that $\frac{y_{f}^{\prime}}{c}$
is odd yielding then $y_{f}$ odd. }{\large\par}

{\large{}For the exceptions of some values of $k$ such that $k\equiv\pm8,\pm16\left(\text{mod}64\right)$,
this procedure does not lead to an odd value of $\frac{y_{f}^{\prime}}{c}$.
For example, for $k=56$, let $k^{\prime}=14$ and $c=2,$yielding
$\left(x_{f},y_{f}^{\prime}\right)=\left(15,4\right)$ and $y_{f}=\frac{y_{f}^{\prime}}{c}=2$.
For $k=72$, let $k^{\prime}=2$ and $c=6,$yielding $\left(x_{f},y_{f}^{\prime}\right)=\left(3,2\right)$.
However, $y_{f}^{\prime}$ is not divisible by $c=6$ and one must
consider not the first fundamental solution of the simple Pell equation
for $k^{\prime}=2$, but the second solution given by (\ref{eq:4})
for $n=2$, yielding $\left(x_{2},y_{2}^{\prime}\right)=\left(17,12\right)$
that gives $y_{f}=\frac{y_{2}^{\prime}}{c}=2$ and finally $\left(x_{f},y_{f}\right)=\left(17,2\right)$.}{\large\par}

{\large{}Furthermore, for some expressions of $k$ in function of
the closest natural square $s^{2}$, one can find general expressions
of $\left(x_{f},y_{f}\right)$ and $\left(X_{f_{i}},Y_{f_{i}}\right)$
in addition to $\left(\pm1,1\right)$ (i.e., for $i>2$) as shown
in Table \ref{tab5: Expressions of k, s, r, }. All these expressions
can easily be demonstrated by replacing the appropriate variables
in (\ref{eq:13}) and (\ref{eq:14}).}{\large\par}

\begin{table}
\caption{\label{tab5: Expressions of k, s, r, }Expressions of $k,s,r,\left(x_{f},y_{f}\right),\left(X_{f_{i}},Y_{f_{i}}\right)$
for $i>2$}

\begin{centering}
\begin{tabular}{|c|c|c|l|}
\hline 
$k$ & $s$ & $r$ & $\left(x_{f},y_{f}\right)\mid\left(X_{f_{i}},Y_{f_{i}}\right)$\tabularnewline
\hline 
\hline 
$s^{2}+1$ & even & 2 & $\left(\pm\left(2s^{2}+1\right),2s\right)\mid\left(-,-\right)$\tabularnewline
\cline{2-4} \cline{3-4} \cline{4-4} 
 & odd & 3 & $\left(\pm\left(2s^{2}+1\right),2s\right)\mid\left(s^{2},s\right)$\tabularnewline
\hline 
$s^{2}+2$ & any & $2^{(a)}$ & $\left(\pm\left(s^{2}+1\right),s\right)\mid\left(-,-\right)$\tabularnewline
\hline 
$s^{2}+4$ & even & $2^{(b)}$ & {\large{}$\left(\pm\left(\frac{s^{2}}{2}+1\right),\frac{s}{2}\right)\mid\left(-,-\right)$}\tabularnewline
\cline{2-4} \cline{3-4} \cline{4-4} 
 & $1\left(\text{mod}4\right)$ & $4^{(c)}$ & {\large{}$\left(\pm\left(\frac{s^{2}\left(s^{2}+3\right)^{2}}{2}+1\right),\frac{s\left(s^{2}+1\right)\left(s^{2}+3\right)}{2}\right)\mid$}\tabularnewline
 &  &  & $\left(\pm\left(s\left(\frac{s^{2}-s+4}{2}\right)-1\right),\frac{s\left(s-1\right)}{2}+1\right)$\tabularnewline
\cline{2-4} \cline{3-4} \cline{4-4} 
 & $3\left(\text{mod}4\right)$ & 4 & {\large{}$\left(\pm\left(\frac{\left(s^{2}+1\right)^{2}\left(s^{2}+4\right)}{2}-1\right),\frac{s\left(s^{2}+1\right)\left(s^{2}+3\right)}{2}\right)\mid$}\tabularnewline
 &  &  & $\left(\pm\left(s\left(\frac{s^{2}+s+4}{2}\right)+1\right),\frac{s\left(s+1\right)}{2}+1\right)$\tabularnewline
\hline 
$s^{2}+8$ & $0\left(\text{mod}4\right)$ & 2 & $\left(\pm\left(\frac{s^{2}}{4}+1\right),\frac{s}{4}\right)\mid\left(-,-\right)$\tabularnewline
\cline{2-4} \cline{3-4} \cline{4-4} 
 & $2\left(\text{mod}4\right)$ & 2 & {\large{}$\left(\pm\left(\frac{s^{2}\left(s^{2}+8\right)}{8}+1\right),s\left(\frac{s^{2}+4}{8}\right)\right)\mid\left(-,-\right)$}\tabularnewline
\hline 
$s^{2}+s$ & any & 2 & $\left(\pm\left(2s+1\right),2\right)\mid\left(-,-\right)$\tabularnewline
\hline 
$s^{2}\widehat{\pm}\frac{2s}{\sigma}$ & $0\left(\text{mod}\sigma\right),\forall\sigma$ odd & $\geq2$ & $\left(\pm\left(\sigma s\widehat{\pm}1\right),\sigma\right)\mid\left(*,*\right)$\tabularnewline
\cline{2-4} \cline{3-4} \cline{4-4} 
 & $0\left(\text{mod}\frac{\sigma}{2}\right),\forall\sigma$ even & $\geq2$ & $\left(\pm\left(\sigma s\widehat{\pm}1\right),\sigma\right)^{(d)}\mid\left(*,*\right)$\tabularnewline
\hline 
$s^{2}+s-1$ & any & $\geq2$ & $\left(*,*\right)\mid\left(\pm\left(2s^{2}+2s-1\right),2s+1\right)^{(e)}$\tabularnewline
\hline 
$s^{2}+s-2$ & $0\left(\text{mod}3\right)$ & >4 & {\large{}$\left(*,*\right)\mid\left(\pm\left(\frac{2s^{2}+4s}{3}-1\right),\frac{2s}{3}+1\right)^{(f)}$}\tabularnewline
\cline{2-4} \cline{3-4} \cline{4-4} 
 & $1\left(\text{mod}3\right)$ & 2 & $\left(*,*\right)\mid\left(-,-\right)$\tabularnewline
\cline{2-4} \cline{3-4} \cline{4-4} 
 & $2\left(\text{mod}3\right)$ & 4 & $\left(*,*\right)\mid\left(\pm\frac{2s^{2}-5}{3},\frac{2\left(s-2\right)}{3}+1\right)$\tabularnewline
\hline 
$s^{2}+s+1$ & $1\left(\text{mod}3\right)$ & 4 & $\left(\pm\left(\frac{2\left(2s+1\right)^{2}}{3}+1\right),4\left(\frac{2\left(s-1\right)}{3}+1\right)\right)\mid$\tabularnewline
 &  &  & $\left(\pm\frac{2s^{2}+2s-1}{3},\frac{2s+1}{3}\right)$\tabularnewline
\cline{2-4} \cline{3-4} \cline{4-4} 
 & $0,2\left(\text{mod}3\right)$ & $\geq4$ & $\left(*,*\right)\mid\left(\pm\left(2s^{2}+2s+1\right),2s+1\right)$\tabularnewline
\hline 
$s^{2}+2s$ & any & 2 & $\left(\pm\left(s+1\right),1\right)\mid\left(-,-\right)$\tabularnewline
\hline 
$s^{2}+2s-1$ & any & 2 & $\left(\pm\left(s^{2}+2s\right),s+1\right)\mid\left(-,-\right)$\tabularnewline
\hline 
$s^{2}+2s-2$ & $2\left(\text{mod}3\right)$ & $2^{(g)}$ & $\left(\pm\frac{2s^{2}+4s-1}{3},\frac{2\left(s+1\right)}{3}\right)\mid\left(-,-\right)$\tabularnewline
\hline 
$s^{2}+2s-3$ & $0\left(\text{mod}4\right)$ & 4 & $\left(\pm\frac{\left(s+1\right)\left(s^{2}+2s-2\right)}{2},\frac{s\left(s+2\right)}{2}\right)\mid\left(\pm\frac{s^{2}+3s-2}{2},\frac{s+2}{2}\right)$\tabularnewline
\cline{2-4} \cline{3-4} \cline{4-4} 
 & $2\left(\text{mod}4\right)$ & $4^{(h)}$ & $\left(\pm\frac{\left(s+1\right)\left(s^{2}+2s-2\right)}{2},\frac{s\left(s+2\right)}{2}\right)\mid\left(\pm\frac{s^{2}+s-4}{2},\frac{s}{2}\right)$\tabularnewline
\cline{2-4} \cline{3-4} \cline{4-4} 
 & odd & $2^{(i)}$ & $\left(\pm\frac{s^{2}+2s-1}{2},\frac{s+1}{2}\right)\mid\left(-,-\right)$\tabularnewline
\hline 
$s^{2}+2s-7$ & $3\left(\text{mod}4\right)$ & $2^{(j)}$ & $\left(\pm\frac{s^{2}+2s-3}{4},\frac{s+1}{4}\right)\mid\left(-,-\right)$\tabularnewline
\cline{2-4} \cline{3-4} \cline{4-4} 
 & $1\left(\text{mod}4\right)$ & 4 & $\left(\pm\frac{s^{2}\left(s^{2}-2\right)+4s\left(s^{2}-3\right)+1}{8},\frac{\left(s^{2}-1\right)\left(s+3\right)}{8}\right)\mid$\tabularnewline
 &  &  & $\left(\pm\frac{2s^{2}+3s-5}{4},\frac{s+1}{2}\right)$\tabularnewline
\hline 
$s^{2}+\frac{\left(3s+1\right)}{2}$ & odd & 2 & $\left(\pm\left(4s+3\right),4\right)\mid\left(-,-\right)$\tabularnewline
\hline 
\end{tabular}
\par\end{centering}
\centering{}$\left(-,-\right)$: no solutions exist as $r=2$;$\left(*,*\right)$:
no apparent pattern; $\widehat{\pm}$: plus/minus sign independent
from other $\pm$ sign; (a) except for $k=51,66$ ($r=3,4$); (b)
except for $k=40$ ($r=4$); (c) except for $k=85$ ($r=8$); (d)
except for $k=\sigma^{2}-1$, with $\sigma$ even ; (e) except for
$k=5,11,55,\ldots$; (f) except for $k=40$; (g) except for $k=78$
($r=4$); (h) except for $k=5$ ($r=2$); (i) except for $k=96$ ($r=4$);
(j) except for $k=136$ ($r=4$)
\end{table}

{\large{}Note that these general expressions for the fundamental solutions
$\left(x_{f},y_{f}\right)$ are valid in all generality for the simple
Pell equation (\ref{eq:14}).}{\large\par}

\subsection{First $r$ solutions of (\ref{eq:1}) for multiple of triangular
numbers}

\label{subsec4.2:First--solutions}

{\large{}Before calculating all solutions of (\ref{eq:1}) yielding
triangular numbers that are multiple of other triangular numbers,
we have to find the first $r$ solutions $\left(\xi_{i},t_{i}\right)$
(with $0\leq i\leq r$) of (\ref{eq:1}), arranged in increasing value
order, i.e., $\xi_{0}=0<\xi_{1}<\ldots<\xi_{i}<\ldots<\xi_{r}$ and
, $t_{0}=0<t_{1}<\ldots<t_{i}<\ldots<t_{r}$, and that correspond
to the $r$ fundamental solutions $\left(X_{f_{i}},Y_{f_{i}}\right)$
of the generalized Pell equation (\ref{eq:13}), with both $X_{f_{i}}$
and $Y_{f_{i}}$ odd or of different parities.}{\large\par}

{\large{}The generalized fundamental solutions $\left(X_{f_{1}},Y_{f_{1}}\right)=\left(1,1\right)$
and $\left(X_{f_{2}},Y_{f_{2}}\right)=\left(-1,1\right)$ provide
respectively, the solutions $\left(\xi_{r},t_{r}\right)$ and $\left(\xi_{r-1},t_{r-1}\right)$
of (\ref{eq:1}) from (\ref{eq:7}), yielding successively
\begin{align*}
\left(X_{1_{1}},Y_{1_{1}}\right) & =\left(X_{f_{1}}x_{f}+kY_{f_{1}}y_{f},X_{f_{1}}y_{f}+Y_{f_{1}}x_{f}\right)\\
 & =\left(x_{f}+ky_{f},y_{f}+x_{f}\right)\\
\left(X_{1_{2}},Y_{1_{2}}\right) & =\left(X_{f_{2}}x_{f}+kY_{f_{2}}y_{f},X_{f_{2}}y_{f}+Y_{f_{2}}x_{f}\right)\\
 & =\left(-x_{f}+ky_{f},-y_{f}+x_{f}\right)
\end{align*}
and
\begin{align}
\left(\xi_{r},t_{r}\right)= & \left(\frac{X_{1_{1}}-1}{2},\frac{Y_{1_{1}}-1}{2}\right)\nonumber \\
= & \left(\frac{x_{f}+ky_{f}-1}{2},\frac{y_{f}+x_{f}-1}{2}\right)\label{eq:18}\\
\left(\xi_{r-1},t_{r-1}\right)= & \left(\frac{X_{1_{2}}-1}{2},\frac{Y_{1_{2}}-1}{2}\right)\nonumber \\
= & \left(\frac{-x_{f}+ky_{f}-1}{2},\frac{-y_{f}+x_{f}-1}{2}\right)\label{eq:19}
\end{align}
Then for $r>2$, the next two generalized fundamental solutions $\left(X_{f_{3}},Y_{f_{3}}\right)$
and $\left(X_{f_{4}},Y_{f_{4}}\right)=\left(-X_{f_{3}},Y_{f_{3}}\right)$
yield respectively $\left(\xi_{1},t_{1}\right)$ and $\left(\xi_{2},t_{2}\right)$.
If both $X_{f_{3}}$ and $Y_{f_{3}}$ are odd, then (\ref{eq:6})
($n=0$) can be used for $\left(\xi_{1},t_{1}\right)$, yielding
\begin{equation}
\left(\xi_{1},t_{1}\right)=\left(\frac{X_{f_{3}}-1}{2},\frac{Y_{f_{3}}-1}{2}\right)\label{eq:20}
\end{equation}
Equation (\ref{eq:6}) could also be used for $\left(\xi_{2},t_{2}\right)$
with $\left(-X_{f_{3}},Y_{f_{3}}\right)$, but it would provide a
negative value for $\xi_{2}$. Instead, we use (\ref{eq:7}) ($n=1$),
giving
\begin{equation}
\left(\xi_{2},t_{2}\right)=\left(\frac{-X_{f_{3}}x_{f}+kY_{f_{3}}y_{f}-1}{2},\frac{-X_{f_{3}}y_{f}+Y_{f_{3}}x_{f}-1}{2}\right)\label{eq:21}
\end{equation}
The next two generalized fundamental solutions $\left(X_{f_{5}},Y_{f_{5}}\right)$
and $\left(X_{f_{6}},Y_{f_{6}}\right)$ $=\left(-X_{f_{5}},Y_{f_{5}}\right)$
yield similarly the next two solutions $\left(\xi_{i},t_{i}\right)$
that are put in the right increasing order.}{\large\par}

{\large{}For example, for $k=13$, $r=4$, $\left(x_{f},y_{f}\right)=\left(649,180\right)$,
$\left(X_{f_{i}},Y_{f_{i}}\right)=\left(\pm1,1\right),\left(\pm25,7\right)$,
(\ref{eq:18}) and (\ref{eq:19}) yield respectively, $\left(\xi_{r},t_{r}\right)=\left(\xi_{4},t_{4}\right)=\left(1494,414\right)$
and $\left(\xi_{r-1},t_{r-1}\right)=\left(\xi_{3},t_{3}\right)=\left(845,234\right)$;
(\ref{eq:20}) and (\ref{eq:21}) yield respectively, $\left(\xi_{1},t_{1}\right)=\left(12,3\right)$
and $\left(\xi_{2},t_{2}\right)=\left(77,21\right)$. }{\large\par}

{\large{}Another example, for $k=46$, $r=6$, $\left(x_{f},y_{f}\right)=(24335,3588)$,
$\left(X_{f_{i}},Y_{f_{i}}\right)=\left(\pm1,1\right),\left(\pm47,7\right),\left(\pm183,27\right)$.
With $\left(X_{f_{1,2}},Y_{f_{1,2}}\right)=\left(\pm1,1\right)$,
(\ref{eq:18}) and (\ref{eq:19}) yield respectively, $\left(\xi_{6},t_{6}\right)=\left(94691,13961\right)$,
$\left(\xi_{5},t_{5}\right)=\left(70356,10373\right)$. With $\left(X_{f_{3,4}},Y_{f_{3,4}}\right)=\left(\pm47,7\right)$,
(\ref{eq:6}) yields $\left(\xi_{1},t_{1}\right)=\left(\frac{X_{f_{3}}-1}{2},\frac{Y_{f_{3}}-1}{2}\right)=\left(23,3\right)$
and (\ref{eq:7}) yields}{\large\par}

{\large{}$\left(\xi_{4},t_{4}\right)=\left(\frac{-X_{f_{3}}x_{f}+kY_{f_{3}}y_{f}-1}{2},\frac{-X_{f_{3}}y_{f}+Y_{f_{3}}x_{f}-1}{2}\right)=\left(5795,854\right)$}{\large\par}

{\large{}Finally, with $\left(X_{f_{5,6}},Y_{f_{5,6}}\right)=\left(\pm183,27\right)$,
(\ref{eq:6}) yields $\left(\xi_{2},t_{2}\right)=\left(\frac{X_{f_{5}}-1}{2},\frac{Y_{f_{5}}-1}{2}\right)=\left(91,13\right)$
and (\ref{eq:7}) yields }{\large\par}

{\large{}$\left(\xi_{3},t_{3}\right)=\left(\frac{-X_{f_{5}}x_{f}+kY_{f_{5}}y_{f}-1}{2},\frac{-X_{f_{5}}y_{f}+Y_{f_{5}}x_{f}-1}{2}\right)=\left(1495,220\right)$.}{\large\par}

{\large{}For the case where $Y_{f_{i}}$ is even, i.e., $k=40,96,120,\ldots$,
(\ref{eq:18}), (\ref{eq:19}) and (\ref{eq:20}) cannot be used with
$\left(X_{f_{1,2}},Y_{f_{1,2}}\right)=\left(\pm1,1\right)$ as both
$k$ and $Y_{f_{i}}$ are even, yielding non-integer solutions for
$\xi$ and $t$. Instead, the other generalized fundamental solution
have to be used with (\ref{eq:7}) ($n=1$) and (\ref{eq:8}) ($n=2$).
For example, for $k=40$, $r=4$, $\left(x_{f},y_{f}\right)=(19,3)$,
$\left(X_{f_{i}},Y_{f_{i}}\right)=\left(\pm1,1\right),\left(\pm11,2\right)$,
(\ref{eq:7}) yields, first, with $\left(X_{f_{3}},Y_{f_{3}}\right)=\left(11,2\right)$,
$\left(X_{1_{3}},Y_{1_{3}}\right)=\left(X_{f_{3}}x_{f}+kY_{f_{3}}y_{f},X_{f_{3}}y_{f}+Y_{f_{3}}x_{f}\right)=\left(449,71\right)$,
yielding $\left(\xi_{2},t_{2}\right)=\left(224,35\right)$, and second,
with $\left(X_{f_{4}},Y_{f_{4}}\right)=\left(-11,2\right)$,}{\large\par}

{\large{}$\left(X_{1_{4}},Y_{1_{4}}\right)=\left(X_{f_{4}}x_{f}+kY_{f_{4}}y_{f},X_{f_{4}}y_{f}+Y_{f_{4}}x_{f}\right)=\left(31,5\right)$,
giving $\left(\xi_{1},t_{1}\right)=\left(15,2\right)$. Next, (\ref{eq:8})
yields, first, with $\left(X_{f_{1}},Y_{f_{1}}\right)=\left(1,1\right)$,}{\large\par}

{\large{}$\left(X_{2_{1}},Y_{2_{1}}\right)=\left(x_{f}^{2}+ky_{f}^{2}+2kx_{f}y_{f},x_{f}^{2}+ky_{f}^{2}+2x_{f}y_{f}\right)=\left(5281,835\right)$,
yielding $\left(\xi_{4},t_{4}\right)=\left(2640,417\right)$, and
second, with $\left(X_{f_{2}},Y_{f_{2}}\right)=\left(-1,1\right)$,
$\left(X_{2_{2}},Y_{2_{2}}\right)=\left(-\left(x_{f}^{2}+ky_{f}^{2}\right)+2kx_{f}y_{f},x_{f}^{2}+ky_{f}^{2}-2x_{f}y_{f}\right)=\left(3839,697\right)$,
yielding $\left(\xi_{3},t_{3}\right)=\left(1919,303\right)$.}{\large\par}

\subsection{All solutions of (\ref{eq:1}) for multiple of triangular numbers}

\label{subsec4.3:All-solutions-of}

{\large{}Once that the first $r$ values of $\left(\xi_{i},t_{i}\right)$
have been found, each corresponding to one of the $r$ generalized
fundamental solutions $\left(X_{f_{i}},Y_{f_{i}}\right)$, the $r$
branches of infinitely many other solutions can be found using either
:}{\large\par}

{\large{}1) the $r$ general solutions (\ref{eq:5}) (assuming a $+$
sign), yielding
\begin{equation}
\xi_{n}+\sqrt{k}t_{n}=\left(\xi_{i}+\sqrt{k}t_{i}+\left(\frac{1+\sqrt{k}}{2}\right)\right)\left(x_{f}+\sqrt{k}y_{f}\right)^{n}-\left(\frac{1+\sqrt{k}}{2}\right)\label{eq:22}
\end{equation}
where $\left(\xi_{i},t_{i}\right)$ must be replaced successively
by the $r$ values of $\left(\xi_{i},t_{i}\right)$; or}{\large\par}

{\large{}2) the first recurrence relation (\ref{eq:9}), yielding
\begin{align}
\left(\xi_{n},t_{n}\right)= & \left(\left(x_{f}\xi_{n-r}+ky_{f}t_{n-r}\right)+\left(\frac{x_{f}+ky_{f}-1}{2}\right),\right.\nonumber \\
 & \left.\left(x_{f}t_{n-r}+y_{f}\xi_{n-r}\right)+\left(\frac{x_{f}+y_{f}-1}{2}\right)\right)\label{eq:23}
\end{align}
where indices of $\xi_{n-r}$ and $t_{n-r}$ (instead of $\xi_{n-1}$
and $t_{n-1}$) in the right part of (\ref{eq:23}) refer to the preceding
values of $\xi$ and $t$ in the same branch of solutions; or,}{\large\par}

{\large{}3) the second recurrence relation (\ref{eq:10}), yielding
\begin{equation}
\left(\xi_{n},t_{n}\right)=\left(2x_{f}\xi_{n-r}-\xi_{n-2r}+\left(x_{f}-1\right),2x_{f}t_{n-r}-t_{n-2r}+\left(x_{f}-1\right)\right)\label{eq:24}
\end{equation}
where indices of $\xi_{n-r},\xi_{n-2r}$ and $t_{n-r},t_{n-2r}$ (instead
of $\xi_{n-1},\xi_{n-2}$ and $t_{n-1},t_{n-2}$) in the right part
of (\ref{eq:24}) refer to the preceding and the one before values
of $\xi$ and $t$ in the same branch of solutions; or,}{\large\par}

{\large{}4) the Chebyshev polynomial solution (\ref{eq:11}), yielding
\begin{align}
\left(\xi_{n},t_{n}\right)= & \left(\left(\xi_{i}+\frac{1}{2}\right)T_{n-1}\left(x_{f}\right)+k\left(t_{i}+\frac{1}{2}\right)y_{f}U_{n-2}\left(x_{f}\right)-\frac{1}{2},\right.\nonumber \\
 & \left.\left(\xi_{i}+\frac{1}{2}\right)y_{f}U_{n-2}\left(x_{f}\right)+\left(t_{i}+\frac{1}{2}\right)T_{n-1}\left(x_{f}\right)-\frac{1}{2}\right)\label{eq:25}
\end{align}
where $\left(\xi_{i},t_{i}\right)$ must be replaced successively
by the $r$ values of $\left(\xi_{i},t_{i}\right)$.}{\large\par}

\subsection{{\large{}Relation between Pell equation solutions and recurrent relations}}

{\large{}\label{subsec4.4:Correspondence-between-Pell}}{\large\par}

{\large{}We can give now a new definition of the rank $r$ introduced
in Section \ref{sec2:Rank-and-recurrent}. The rank $r$ is the number
of fundamental solutions $\left(X_{f_{i}},Y_{f_{i}}\right)$ of the
generalized Pell equation (\ref{eq:13}), with $X_{f_{i}}$ odd and
$Y_{f_{i}}$ odd or even (if $y_{f}$ is not even) , with $r\leq\rho$,
the total number of generalized solution of (\ref{eq:13}).}{\large\par}

{\large{}Furthermore, we see that the second recurrent relations (\ref{eq:24})
for both $\xi_{n}$ and $t_{n}$ have $x_{f}$ as the only parameter,
and that the two relations are independent from the value of $k$
and $y_{f}$. This fundamental solution $x_{f}$ of the simple Pell
equation (\ref{eq:14}) acts like a constant of the problem for each
value of $k$. Note further that summing the expressions of $t_{r}$
and $t_{r-1}$ in (\ref{eq:18}) and (\ref{eq:19}) yields $t_{r}+t_{r-1}=x_{f}-1$.
As this sum $t_{r}+t_{r-1}$ was already defined in (\ref{eq:3.2}),
the constant $\kappa$ is related to $x_{f}$
\begin{equation}
\kappa=x_{f}-1\label{eq:34}
\end{equation}
Furthermore, (\ref{eq:34}) yields also that $y_{f}$ is related to
the difference $\delta=t_{r}-t_{r-1}$ through the simple Pell equation
\ref{eq:14} $\left(\kappa+1\right)^{2}-ky_{f}^{2}=1$, which is verified
if $y_{f}^{2}=\kappa^{2}-4t_{r}t_{r-1}=\left(t_{r}-t_{r-1}\right)^{2}=\delta^{2}$,
giving
\begin{equation}
\delta=y_{f}\label{eq:35}
\end{equation}
for all non-square values of $k$ (except for some values of $k$
such that $k\equiv0\left(\text{mod}8\right)$, see further). Replacing
in the simple Pell equation $\left(\kappa+1\right)^{2}-k\delta^{2}=1$
yields the condition between the sum and the difference of $t_{r}$
and $t_{r-1}$
\begin{equation}
\delta=\sqrt{\frac{\kappa^{2}+\kappa}{k}}\label{eq:36}
\end{equation}
With the exception of $k=56$, 72, 112, 184, 240, 248, 264, 272, 376,
..., i.e., for some values of $k$ such that $k\equiv\pm8,\pm16\left(\text{mod}64\right)$
(for which (\ref{eq:35}) is valid), the relation (\ref{eq:35}) is
not valid for the other cases of $k\equiv0\left(\text{mod}8\right)$.
In these cases, $\delta>y_{f}$ and one must find the next pair of
solutions to the simple Pell equation by (\ref{eq:4}) for $n=2$,
i.e., $\left(x_{2},y_{2}\right)=\left(x_{f}^{2}+ky_{f}^{2},2x_{f}y_{f}\right)$.
Then for these cases,
\begin{align}
\kappa & =x_{f}^{2}+ky_{f}^{2}-1\label{eq:37}\\
\delta & =2x_{f}y_{f}\label{eq:38}
\end{align}
Finally, replacing $x_{f}$ in (\ref{eq:24}) from \ref{eq:34} yields
\begin{equation}
\left(\xi_{n},t_{n}\right)=\left(2\left(\kappa+1\right)\xi_{n-r}-\xi_{n-2r}+\kappa,2\left(\kappa+1\right)t_{n-r}-t_{n-2r}+\kappa\right)\label{eq:26}
\end{equation}
which are the same recurrent relations given in (\ref{eq:3.3}) and
(\ref{eq:3.3-1}).}{\large\par}

\section{{\large{}Conclusions}}

\label{sec5:Conclusions}

{\large{}We have shown that the problem of finding all triangular
numbers that are multiples of other triangular numbers with non-square
integer multiplier $k$ can be solved using solutions of Pell equations
with a simple change of variables. Only those $r$ fundamental solutions
$\left(X_{f_{i}},Y_{f_{i}}\right)$ of the generalized Pell equation
with $X_{f_{i}}$ odd and $Y_{f_{i}}$ odd or even (if $y_{f}$ is
not even) provide solutions to the problem of finding triangular numbers
that are multiple of other triangular numbers. General expressions
of fundamental solutions of the Pell equations are given for some
values of the multiplier $k$ in function of the closest natural square
values $s^{2}$. Many infinitely solutions are then found on $r$
branches corresponding to each of the $r$ generalized fundamental
solutions $\left(X_{f_{i}},Y_{f_{i}}\right)$ and these solutions
can be found either by a general relation involving $\sqrt{k}$, or
by a first set of recurrent relations, or by a second set of recurrent
relations, or by Chebyshev polynomial solutions. Among these, the
second set of recurrent relations are found to be the same as those
found previously without using the Pell equation solving method. }{\large\par}

{\large{}Furthermore, the number $r$ of generalized fundamental solutions
$\left(X_{f_{i}},Y_{f_{i}}\right)$ with $X_{f_{i}}$ odd and $Y_{f_{i}}$
odd or even (if $y_{f}$ is not even) corresponds to the rank of these
second set recurrent relations. Finally, the two constants $\kappa=t_{r}+t_{r-1}$
and $\delta=t_{r}-t_{r-1}$ are also related to respectively the fundamental
solutions $x_{f}$ and $y_{f}$ of the simple Pell equation, as $\kappa=x_{f}-1$
and $\delta=y_{f}$ or $\delta=2x_{f}y_{f}$. These are an unexpected
result as this means that the fundamental solutions of the simple
Pell equation, in all its generality, are related to constants in
recurrent relations of the problem of finding triangular numbers multiple
of other triangular numbers.}{\large\par}

\end{document}